%% file: nm_rom_neurips.tex
\documentclass{article}

     \PassOptionsToPackage{numbers, compress}{natbib}

    \usepackage[final]{ml4eng_2020}

\usepackage[utf8]{inputenc} 
\usepackage[T1]{fontenc}    
\usepackage{hyperref}       
\usepackage{url}            
\usepackage{booktabs}       
\usepackage{amsfonts}       
\usepackage{nicefrac}       
\usepackage{microtype}      

\usepackage{upgreek}
\usepackage[mathscr]{eucal}
\usepackage{amsmath,amssymb,amsopn,mathtools}
\usepackage{graphicx}
\usepackage{color}
\usepackage{relsize}
\usepackage{enumitem}
\usepackage[normalem]{ulem}
\usepackage{verbatim}
\usepackage{footmisc}
\usepackage{soul}
\usepackage{subfigure, algorithmic, algorithm}
\usepackage{multirow, booktabs}
\usepackage{siunitx}
\usepackage{xr}

\hypersetup{
	colorlinks=true,   
	citecolor=blue,     
	filecolor=blue,     
	linkcolor=red,    
	urlcolor=blue}      
	
\definecolor{Blue}{rgb}{0,0,1}
\definecolor{Red}{rgb}{1,0,0}
\definecolor{Green}{rgb}{0,1,0}
\definecolor{Bronze}{rgb}{0.8,0.5,0.2}
\definecolor{Violet}{rgb}{0.54,0.17,0.89}

\input{commands}

\title{Efficient nonlinear manifold reduced order model}

\author{
  Youngkyu Kim \\
  Mechanical Engineering\\
  University of California, Berkeley\\
  Berkeley, CA 94720 \\
  \texttt{youngkyu$\_$kim@berkeley.edu} \\
  \And
  Youngsoo Choi \\
  Center for Applied Scientific Computing \\
  Lawrence Livermore National Laboratory \\
  Livermore, CA 94550 \\
  \texttt{choi15@llnl.gov} \\
  \AND
  David Widemann \\
  Computational Engineering Division \\
  Lawrence Livermore National Laboratory \\
  Livermore, CA 94550 \\
  \texttt{widemann1@llnl.gov} \\
  \And
  Tarek Zohdi \\
  Mechanical Engineering\\
  University of California, Berkeley\\
  Berkeley, CA 94720 \\
  \texttt{zohdi@berkeley.edu} \\
}

\begin{document}

\maketitle

\begin{abstract}
Traditional linear subspace reduced order models (LS-ROMs) are able to
accelerate physical simulations, in which the intrinsic solution space falls
into a subspace with a small dimension, i.e., the solution space has a small
Kolmogorov $n$-width. However, for physical phenomena not of this type, such
as advection-dominated flow phenomena, a low-dimensional linear subspace
poorly approximates the solution. To address cases such as these, we have
developed an efficient nonlinear manifold ROM (NM-ROM), which can better
approximate high-fidelity model solutions with a smaller latent space
dimension than the LS-ROMs. Our method takes advantage of the existing
numerical methods that are used to solve the corresponding full order models
(FOMs).  The efficiency is achieved by developing a hyper-reduction technique
in the context of the NM-ROM.  Numerical results show that neural networks can
learn a more efficient latent space representation on advection-dominated data
from 2D Burgers' equations with a high Reynolds number.  A speed-up of up to
$11.7$ for 2D Burgers' equations is achieved with an appropriate treatment of
the nonlinear terms through a hyper-reduction technique.  
\end{abstract}

\section{Introduction}\label{sec:intro}
Physical simulations are influencing developments in science, engineering, and
technology more rapidly than ever before. However, high-fidelity, forward
physical simulations are computationally expensive and, thus, make intractable
many decision-making applications, such as design optimization, inverse
problems, optimal controls, and uncertainty quantification. These applications
require many forward simulations to explore the parameter space in the outer
loop. To compensate for the computational expense issue, many surrogate models
have been developed: from simply using interpolation schemes for specific
quantity of interests to physics-informed surrogate models. This paper focuses
on the latter because a physics-informed surrogate model is more robust in
predicting physical solutions than the simple interpolation schemes. 

Among many types of physics-informed surrogate models, the projection-based
linear subspace reduced order models (LS-ROMs) take advantage of both the known
governing equation and data with linear subspace solution representation
\cite{benner2015survey}. Although LS-ROMs have been successfully applied to many
forward physical problems \cite{ghasemi2015localized, jiang2019implementation,
yang2016fast, yang2017efficient, zhao2014pod, cstefuanescu2013pod,
mordhorst2017pod, dimitriu2013application, antil2012reduced} and partial
differential equation(PDE)-constrained optimization problems
\cite{amsallem2015design, choi2020gradient, choi2019accelerating, fu2018pod},
the linear subspace solution representation suffers from not being able to
represent certain physical simulation solutions with a small basis dimension,
such as advection-dominated or sharp gradient solutions. This is because LS-ROMs
work only for physical problems, in which the intrinsic solution space falls
into a subspace with a small dimension, i.e., the solution space has a small
Kolmogorov $n$-width. Although there have been many attempts to resolve these
shortcomings of LS-ROMs with various methods, \cite{abgrall2016robust,
reiss2018shifted, carlberg2013gnat, carlberg2018conservative, choi2020sns,
xiao2014non, constantine2012reduced, choi2020space, choi2019space,
taddei2020space, carlberg2015adaptive, rim2018transport, parish2019windowed,
peherstorfer2018model, welper2020transformed}, all these approaches are still
based on the linear subspace solution representation. We transition to a
nonlinear, low-dimensional manifold to approximate the solution better than
linear methods.

There are many works available in the current literature that looked into the
nonlinear manifold solution represenation, using neural networks (NNs) as
surrogates for physical simulations \cite{lagaris1998artificial,
dissanayake1994neural, van1995neural, meade1994numerical, raissi2019physics,
chen2018neural, khoo2017solving, long2018pde, beck2019machine, zhu2019physics,
berg2018unified, sirignano2018dgm, han2018solving, lu2019deeponet,
lu2019deepxde, pang2019fpinns, zhang2019quantifying, weinan2018deep,
he2018relu}. However, these methods do not take advantage of the existing
numerical methods for high-fidelity physical simulations. Recently, a neural
network-based ROM is developed in \cite{lee2020model}, where the weights and
biases are determined in the training phase and the existing numerical methods
are utilized in their models. The same technique is extended to preserve the
conserved quantities in the physical conservation laws \cite{lee2019deep}.
However, their approaches do not achieve any speed-up because the nonlinear
terms that still scale with the corresponding FOM size need to be updated
every time step or Newton step.                           

We present a fast and accurate physics-informed neural network ROM with a
nonlinear manifold solution representation, i.e., the nonlinear manifold ROM
(NM-ROM). We train a shallow masked autoencoder with solution data from the
corresponding FOM simulations and use the decoder as the nonlinear manifold
solution representation. Our NM-ROM is different from the aformentioned
physics-informed neural networks in that we take advantage of the existing
numerical methods of solving PDE in our approach and a considerable speed-up is
achieved.

\section{Full order model}\label{sec:FOM}
A parameterized nonlinear dynamical system is considered, characterized by
a system of nonlinear ordinary differential equations (ODEs), which can be
considered as a resultant system from semi-discretization of Partial
Differential Equations (PDEs) in space domains
 \begin{equation} \label{eq:fom}
  \frac{d\sol}{dt} = \flux(\sol,t; \param),\quad\quad
  \sol(0;\param) = \solArg{0}(\param),
 \end{equation}     

where $t\in \timeDomain$ denotes time with the final time
$\totaltime\in\RRplus{}$, and $\sol(t;\param)$ denotes the time-dependent,
parameterized state implicitly defined as the solution to
problem~\eqref{eq:fom} with $\sol:\timeDomain\times \paramDomain\rightarrow
\RR{\nspacedof}$.  Further, $\flux: \RR{\nspacedof} \times \timeDomain
\times \paramDomain \rightarrow \RR{\nspacedof}$ with
$(\solDummy,\timeDummy;\paramDummy)\mapsto\flux(\solDummy,\timeDummy;\paramDummy)
$ denotes the time derivative of $\sol$, which we assume to be nonlinear
in at least its first argument.  The initial state is denoted by
$\solArg{0}:\paramDomain\rightarrow \RR{\nspacedof}$, and $\param \in
\paramDomain$ denotes parameters in the domain
$\paramDomain\subseteq\RR{\nparam}$. 

A uniform time discretization is assumed throughout the paper, characterized
by time step $\dt\in\RRplus{}$ and time instances $\timeArg{n} =
\timeArg{n-1} + \dt$ for $n\in\nat{\ntimedof}$ with $\timeArg{0} = 0$,
$\ntimedof\in\natNo$, and $\nat{N}\defeq\{1,\ldots,N\}$.  To avoid
notational clutter, we introduce the following time discretization-related
notations: $\solArg{n} \defeq \solFuncArg{n}$, $\solapproxArg{n} \defeq
\solapproxFuncArg{n}$, $\redsolapproxArg{n} \defeq \redsolapproxFuncArg{n}$,
and $\fluxArg{n} \defeq \flux(\solFuncArg{n},t^{n}; \param)$, where $\sol$, $\solapprox$, $\redsolapprox$ and $\flux$ are defined in. 
The implicit backward Euler (BE)\footnote{Other time integrators can be used in
our NM-ROMs.} time integrator numerically solves Eq.~\eqref{eq:fom}, by solving
the following nonlinear system of equations, i.e.,  
  $\solArg{n} - \solArg{n-1} = \dt\fluxArg{n}$,
for $\solArg{n}$ at $n$-th time step.  The corresponding residual function is
defined as
\begin{align}\label{eq:residual_BE} 
\begin{split}
  \resn_{\BE}(\solArg{n};\solArg{n-1},\param) &\defeq 
  \solArg{n} - \solArg{n-1} -\dt\fluxArg{n}.
\end{split}
\end{align}
\section{Nonlinear manifold reduced order model (NM-ROM)}\label{sec:NM-ROM}
The NM-ROM applies solution representation using a nonlinear manifold
$\spatialSubspace \defeq \{\scaledDecoder\left(\reddummy\right)|\reddummy \in
\RR{\nbasisspace}\}$, where $\scaledDecoder: \RR{\nbasisspace} \rightarrow
\RR{\nspacedof}$ with $\nbasisspace\ll\nspacedof$ denotes a nonlinear function
that maps a latent space of dimension $\nbasisspace$ to the full order model
space of dimension, $\nspacedof$. That is, the NM-ROM approximates the solution
in a trial manifold as 
\begin{equation}\label{eq:spatialNMROMsolution} 
  \sol \approx \solapprox= \solArg{ref} + \scaledDecoder \left(\redsolapprox
  \right).
\end{equation} 
The construction of the nonlinear function, $\scaledDecoder$, is explained in
Section~\ref{sec:NN}.  By plugging Eq.~\eqref{eq:spatialNMROMsolution} into
Eq.~\eqref{eq:residual_BE}, the residual function at $n$th time step becomes
\begin{align}\label{eq:trialManifold_residual_BE} 
\begin{split}
 \resRedApproxArg{n}_{\BE}(\redsolapproxArg{n};\redsolapproxArg{n-1},\param)
  &\defeq
  \resn_{\BE}(\solArg{ref}+\scaledDecoder\left(\redsolapproxArg{n}\right);
  \solArg{ref}+\scaledDecoder\left(\redsolapproxArg{n-1}\right),\param) \\ &=
  \scaledDecoder\left(\redsolapproxArg{n}\right)-\scaledDecoder\left(\redsolapproxArg{n-1}\right)
  -\dt\flux(\solArg{ref}+\scaledDecoder\left(\redsolapproxArg{n}\right),t_n;\param),
\end{split}
\end{align}
which is an over-determined system that we close with the least-squares
Petrov--Galerkin (LSPG) projection. That is, we minimize the squared norm of the
residual vector function at every time step:
\begin{align} \label{eq:manifoldOpt1}
  \begin{split} 
    \redsolapproxArg{n} = \argmin{\reddummy\in\RR{\nbasisspace}} \quad&
    \frac{1}{2} \left \|\resRedApproxArg{n}_{\BE}(\reddummy;\redsolapproxArg{n-1},\param)
    \right \|_2^2.
  \end{split} 
\end{align}
The Gauss--Newton method with the starting point $\redsolapproxArg{n-1}$ is
applied to solve the minimization problem~\eqref{eq:manifoldOpt1}. However, 
the nonlinear residual vector, $\resRedApproxArg{n}_{\BE}$, scales with FOM size
and it needs to be updated every time the argument of the function changes,
which occurs either at every time step or Gauss--Newton step.  More
specifically, if the backward Euler time integrator is used,
$\scaledDecoder\left(\redsolapproxArg{n}\right)$,
$\flux(\solArg{ref}+\scaledDecoder\left(\redsolapproxArg{n}\right),t;\param)$,
and their Jacobians need to be updated whenever $\redsolapproxArg{n}$
changes.  Without any special
treatment on the nonlinear residual term, no speed-up can be expected. Thus, we
apply a hyper-reduction to eliminate the scale with FOM size in the nonlinear
term evaluations (see Section~\ref{sec:HR}). Finally, we denote this
non-hyper-reduced NM-ROM as NM-LSPG.
\section{Shallow masked autoencoder}\label{sec:NN}
The nonlinear function, $\scaledDecoder$, is the decoder
$\unscaledDecoder$ of an autoencoder in the form of a feedforward neural
network. The autoencoder compresses FOM solutions of Eq.~\eqref{eq:fom} with an
encoder $\unscaledEncoder$ and decompresses back to reconstructed FOM solution
with an decoder $\unscaledDecoder$.  The autoencoder is trained to reconstruct
the FOM solutions of Eq.~\eqref{eq:fom} by minimizing the mean square error
between original and reconstructed FOM solutions.  Therefore, the dimension of
the encoder input and the decoder output is $\nspacedof$ and the dimension of
the encoder $\unscaledEncoder$ output and the decoder $\unscaledDecoder$ input
is $\nbasisspace$. 

\begin{figure}[!htbp]                                                           
  \centering                                                                    
  \subfigure[Without masking]{                                                  
  \includegraphics[width=0.45\textwidth]{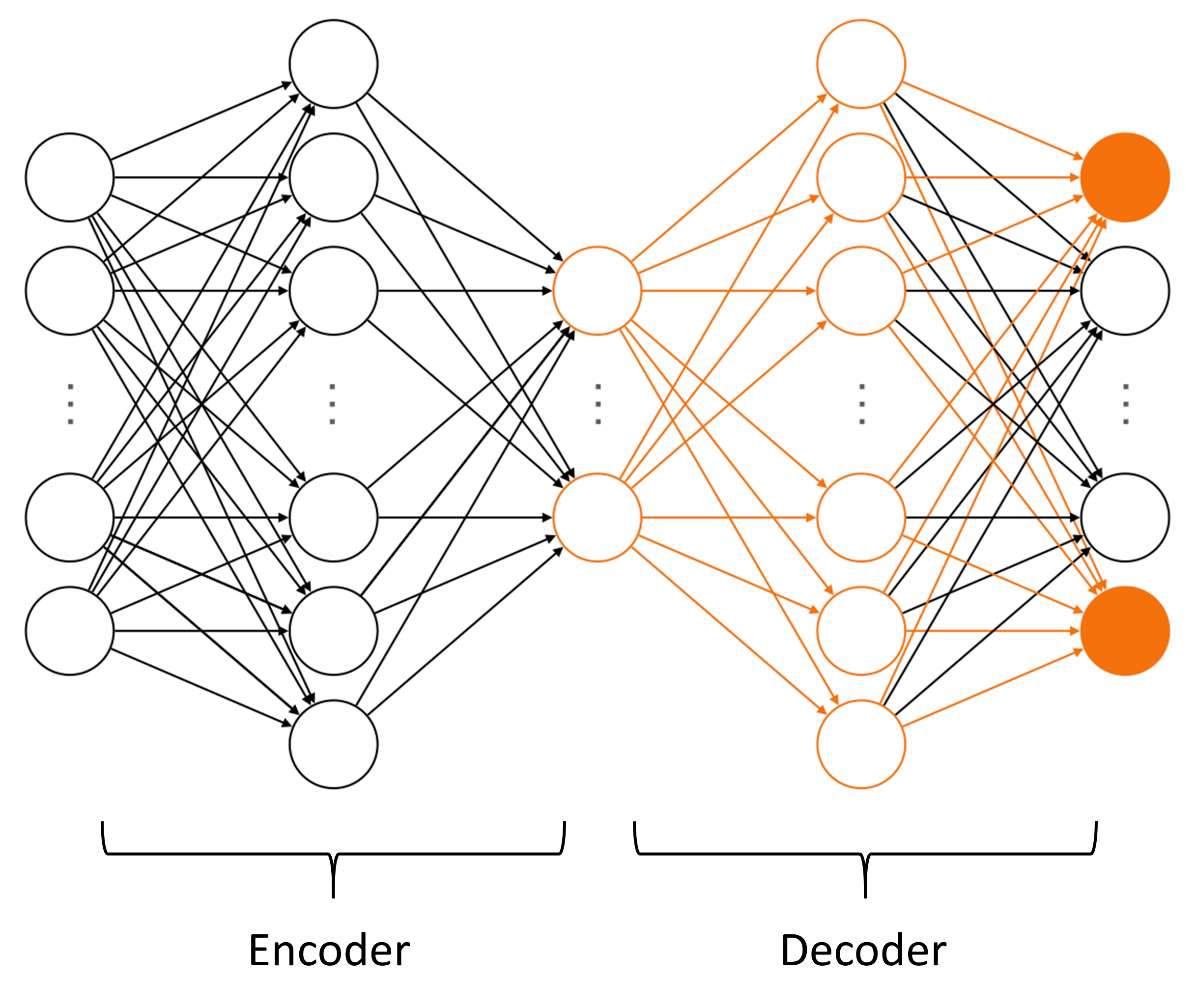}                 
  }                                                                             
  \subfigure[With masking]{                                                     
  \includegraphics[width=0.45\textwidth]{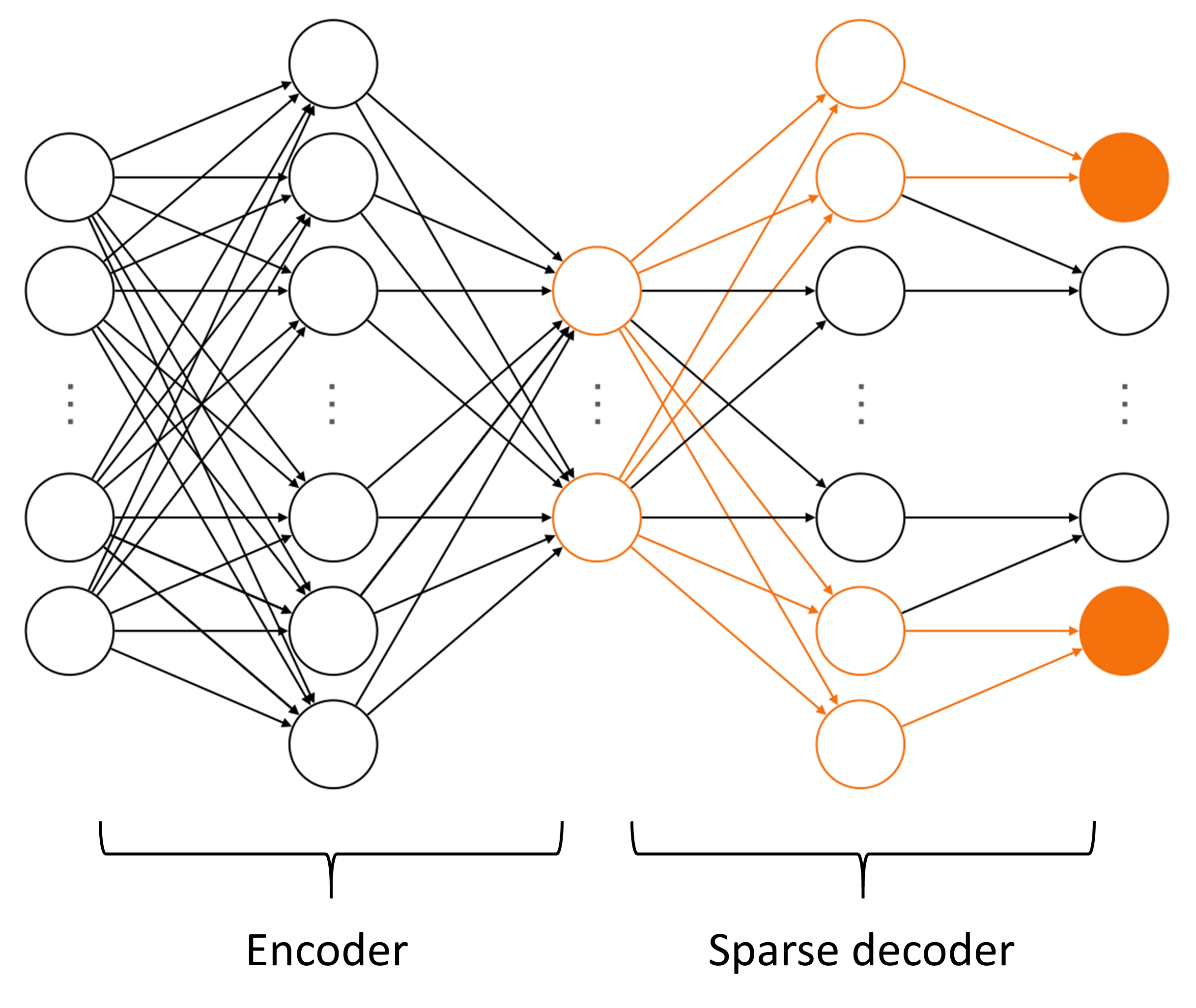}                
  }                                                                             
  \caption{Three-layer encoder/decoder architectures: (a) unmasked and (b) after
  the sparsity mask is applied. Nodes and edges in orange color represent
  active path in the subnet that stems from the sampled outputs that are marked
  as the orange disks. Note that the masked shallow neural network has a sparser
  structure than the unmasked one. A sparser structure leads to a more
  efficient model.} 
  \label{fg:threeLayerAE}                                                       
\end{figure}

We intentionally use a non-deep neural network, i.e., three-layer autoencoder,
for the decoder to achieve an efficiency that is requried by the hyper-reduction
(see Section~\ref{sec:HR}). More specifically, the first layers of the encoder
$\unscaledEncoder$ and decoder $\unscaledDecoder$ are fully-connected layers,
where the nonlinear activation functions are applied and the last layer of the
encoder $\unscaledEncoder$ is fully-connected layer with no activation
functions. The last layer of the decoder $\unscaledDecoder$ is
sparsely-connected layer with no activation functions. The sparsity is
determined by a mask matrix.  These network architectures are shown in
Fig.~\ref{fg:threeLayerAE}.

There is no way to determine hidden layer sizes \textit{a priori}. If the
number of learnable parameters is not enough, the decoder is not able to
represent the nonlinear manifold well.  On the other hand, too many learnable
parameters may result in over-fitting, so the decoder is not able to generalize
well, which means the trained decoder cannot be used for problems whose data is
unseen, i.e., the predictive case.
To avoid over-fitting, we first divide the training
data into train and validation datasets.  Then, the autoencoder is trained
using the train dataset and tested for the generalization using the
validation dataset.  If the mean squared errors on the validation and train datasets are very
different, then over-fitting has occurred. We then reduce the size of the hidden
layer and re-train the model \cite{kramer1991nonlinear}.

\section{Hyper-reduction}\label{sec:HR}
The hyper-reduction techniques are developed to eliminate the FOM scale
dependecy in nonlinear terms \cite{chaturantabut2010nonlinear,
drmac2016new, drmac2018discrete, carlberg2013gnat, choi2020sns}, which is
essential to acheive an efficiency in our NM-ROM. We follow the gappy POD
approach \cite{everson1995karhunen}, in which the nonlinear residual term is
approximated as
 \begin{equation}\label{eq:ResApprox}
   \resRedApproxArg{} \approx \basismatres\redres,
 \end{equation}
where $\basismatres \defeq
[\basisresvecArg{1},\ldots,\basisresvecArg{\nbasisres} ] \in
\RR{\ndof\times\nbasisres}$, $\nbasisspace \leq \nbasisres \ll \ndof$, denotes
the residual basis matrix and $\redres \in \RR{\nbasisres}$ denotes the
generalized coordinates of the nonlinear residual term.  Here,
$\resRedApproxArg{}$ represents a residual vector function, e.g., the backward
Euler residual, $\resRedApproxArg{n}_{BE}$, defined in
Eq.~\eqref{eq:trialManifold_residual_BE}.  We use the singular value
decomposition of the FOM solution snapshot matrix to construct $\basismatres$,
which is justified in \cite{choi2020sns}.  In order to find $\redres$, we apply
a sampling matrix $\samplemat \defeq [\unitvecArg{p_1}, \ldots,
\unitvecArg{p_{\nressample}}]^T \in \RR{\nressample\times\ndof}$, $\nbasisspace
\leq \nbasisres \leq \nressample \ll \ndof$ on both sides of
\eqref{eq:ResApprox}. The vector,  $\unitvecArg{p_i}$, is the $p_i$th column of
the identity matrix $\identity{\ndof}\in\RR{\ndof\times\ndof}$.  Then the
following least-squares problem is solved: 
\begin{align} \label{eq:ResApprox_least-squares}
  \begin{split} 
    \redres := \argmin{\reddummy\in\RR{\nbasisres}} \quad&
    \frac{1}{2} \left \|\samplemat(\resRedApproxArg{} - \basismatres\reddummy)
    \right \|_2^2.
  \end{split} 
\end{align}
The solution to Eq.~\eqref{eq:ResApprox_least-squares} is given as $ \redres =
(\samplemat\basismatres)^\dagger\samplemat\resRedApproxArg{}$, where the
Moore--Penrose inverse of a matrix $\weightmat \in \RR{\nressample \times
\nbasisres}$ with full column rank is defined as $\weightmat^{\dagger} :=
(\weightmat^T\weightmat)^{-1}\weightmat^T$. Therefore, Eq.~\eqref{eq:ResApprox}
becomes $\resRedApproxArg{} \approx \obliqueprojector \resRedApproxArg{}$,
where $\obliqueprojector:= \basismatres
(\samplemat\basismatres)^\dagger\samplemat$ is the oblique projection matrix.
We do not construct the sampling matrix $\samplematNT$.  {\it Instead, it
maintains the sampling indices $\{p_1,\ldots,p_{\nbasisflux}\}$ and
corresponding rows of $\basismatres$ and $\resRedApproxArg{}$.} This enables
hyper-reduced ROMs to achieve a speed-up. 
The sampling indices (i.e., $\samplematNT$) can be determined by Algorithm 3 of
\cite{carlberg2013gnat} for computational fluid dynamics problems and Algorithm
5 of \cite{carlberg2011efficient} for other problems. 

The hyper-reduced residual,  $\obliqueprojector\resRedApproxArg{n}_{BE}$, is
used in the minimization problem in Eq.~\eqref{eq:manifoldOpt1}: 
  \begin{align} \label{eq:NM-LSPG-HR-spOpt1}
    \begin{split} 
      \redsolapproxArg{n} = \argmin{\reddummy\in\RR{\nbasisspace}} \quad&
      \frac{1}{2} \left \|\ (\samplemat\basismatres)^\dagger\samplemat
      \resRedApproxArg{n}_{\BE}(\reddummy;\redsolapproxArg{n-1},\param)
      \right \|_2^2.
    \end{split} 
  \end{align}
Note that the pseudo-inverse $(\samplemat\basismatres)^\dagger$ can be
pre-computed.  Due to the definition of $\resRedApproxArg{n}_{\BE}$
in Eq.~\eqref{eq:trialManifold_residual_BE}, the sampling matrix $\samplematNT$
needs to be applied to the following two terms:
$\scaledDecoder\left(\redsolapproxArg{n}\right) -
\scaledDecoder\left(\redsolapproxArg{n-1}\right)$ and
$\flux(\solArg{ref}+\scaledDecoder\left(\redsolapproxArg{n}\right),t;\param)$ at
every time step. The first term,
$\samplemat(\scaledDecoder\left(\redsolapproxArg{n}\right) -
\scaledDecoder\left(\redsolapproxArg{n-1}\right))$, requires that only selected
outputs of the decoder be computed. Furthermore, for the second term, the nonlinear residual elements that are selected by the sampling matrix need
to be computed. This implies that we have to keep track of the outputs of
$\scaledDecoder$ that are needed to compute the selected nonlinear residual
elements by the sampling matrix, which is usually a larger set than the outputs
that are selected solely by the sampling matrix. Therefore, we
build a subnet that computes only the outputs of the decoder that are
required to compute the nonlinear residual elements. Such outputs are demonstrated as the solid oragne disks and the corresponding subnet is depicted in Fig.~\ref{fg:threeLayerAE}(b). Finally, we denote
this hyper-reduced NM-ROM as NM-LSPG-HR. 

\section{2D Burgers' equation}\label{sec:2dburgers}
\begin{figure}[!htbp]
  \centering
  \subfigure[FOM]{
  \includegraphics[width=0.25\textwidth]{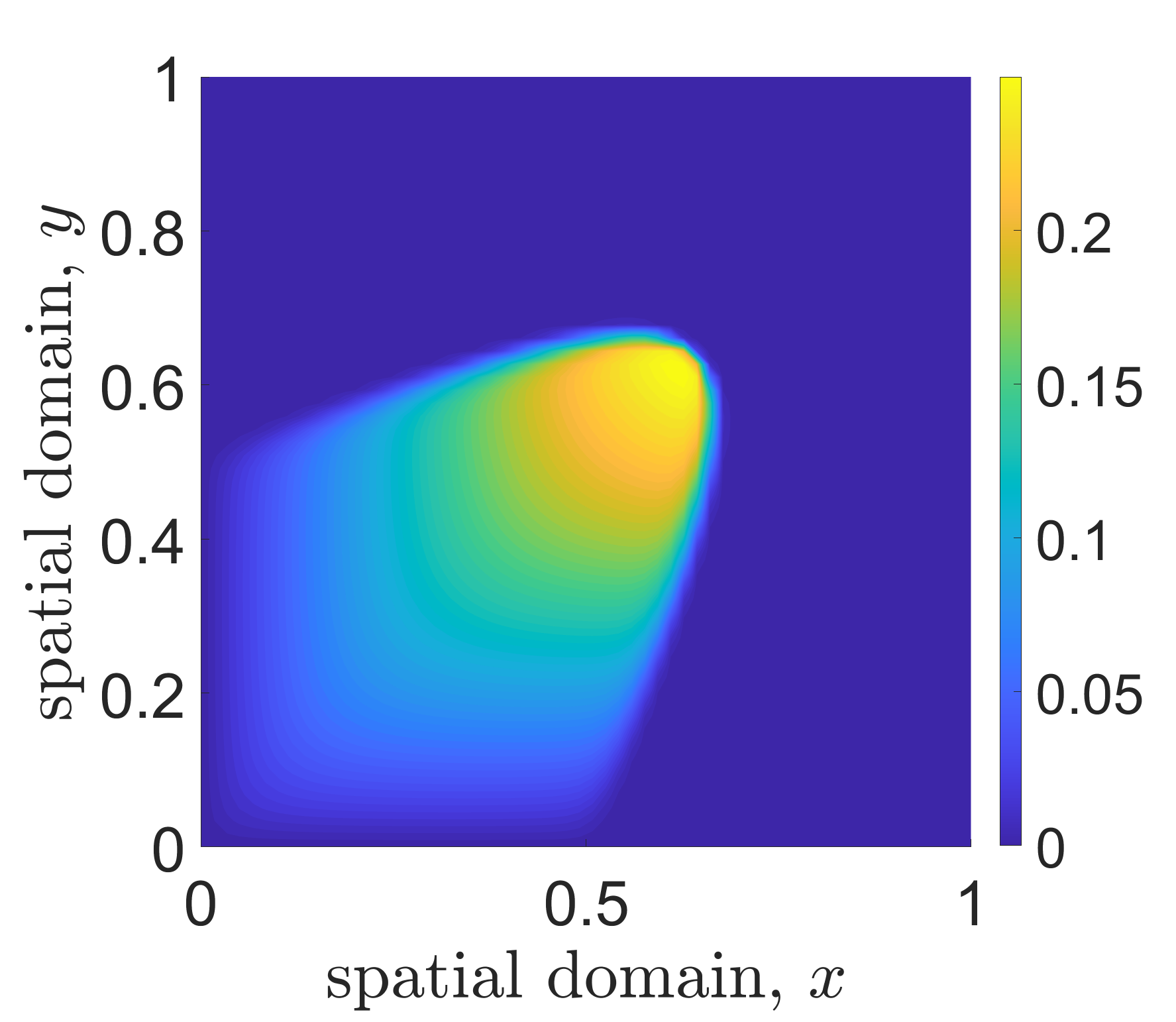}
}~~~~~~~ \hspace{-22pt}
\subfigure[NM-LSPG-HR]{
  \includegraphics[width=0.25\textwidth]{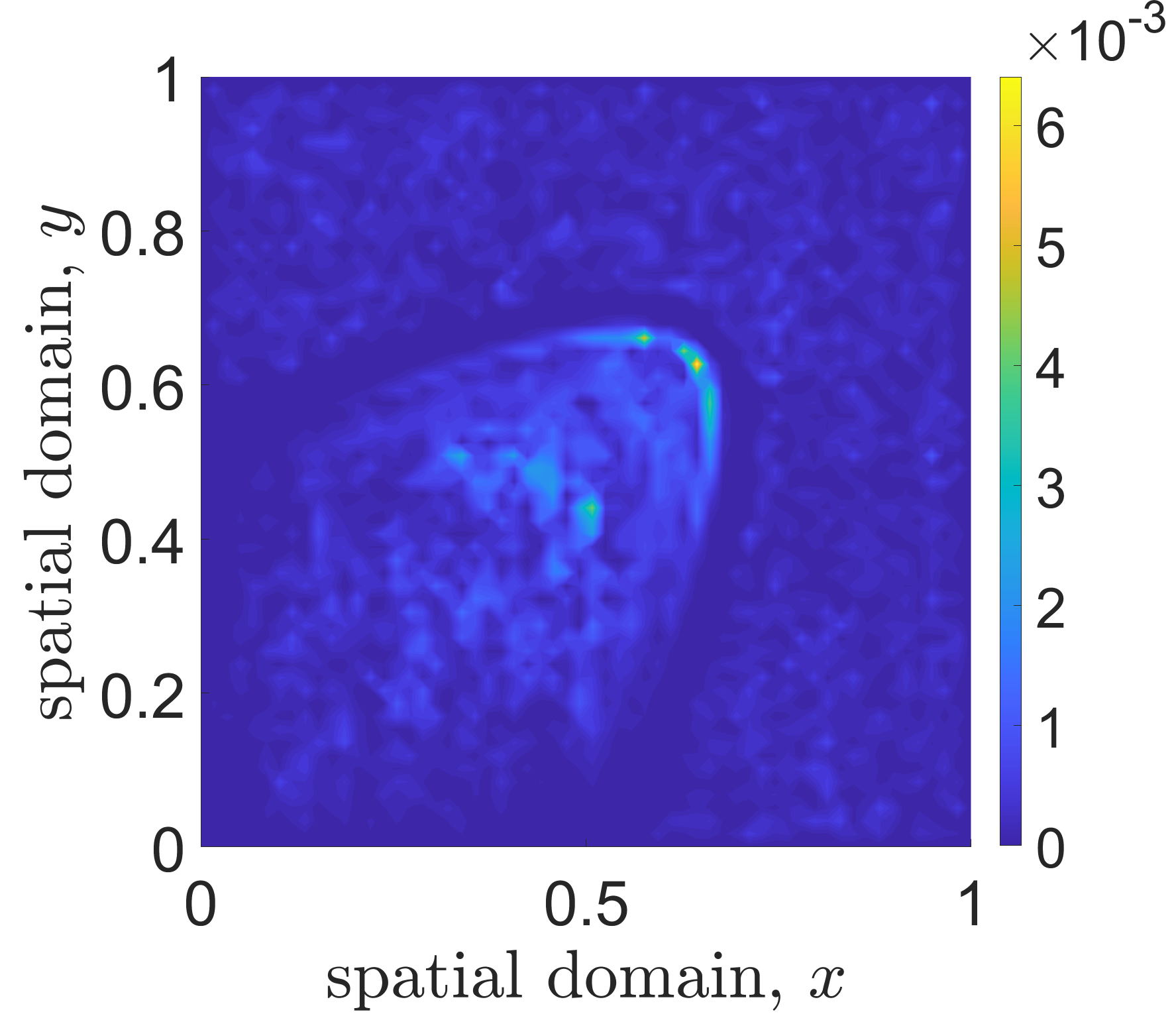}
}~~~~~~~ \hspace{-22pt}
\subfigure[LS-LSPG-HR]{
  \includegraphics[width=0.25\textwidth]{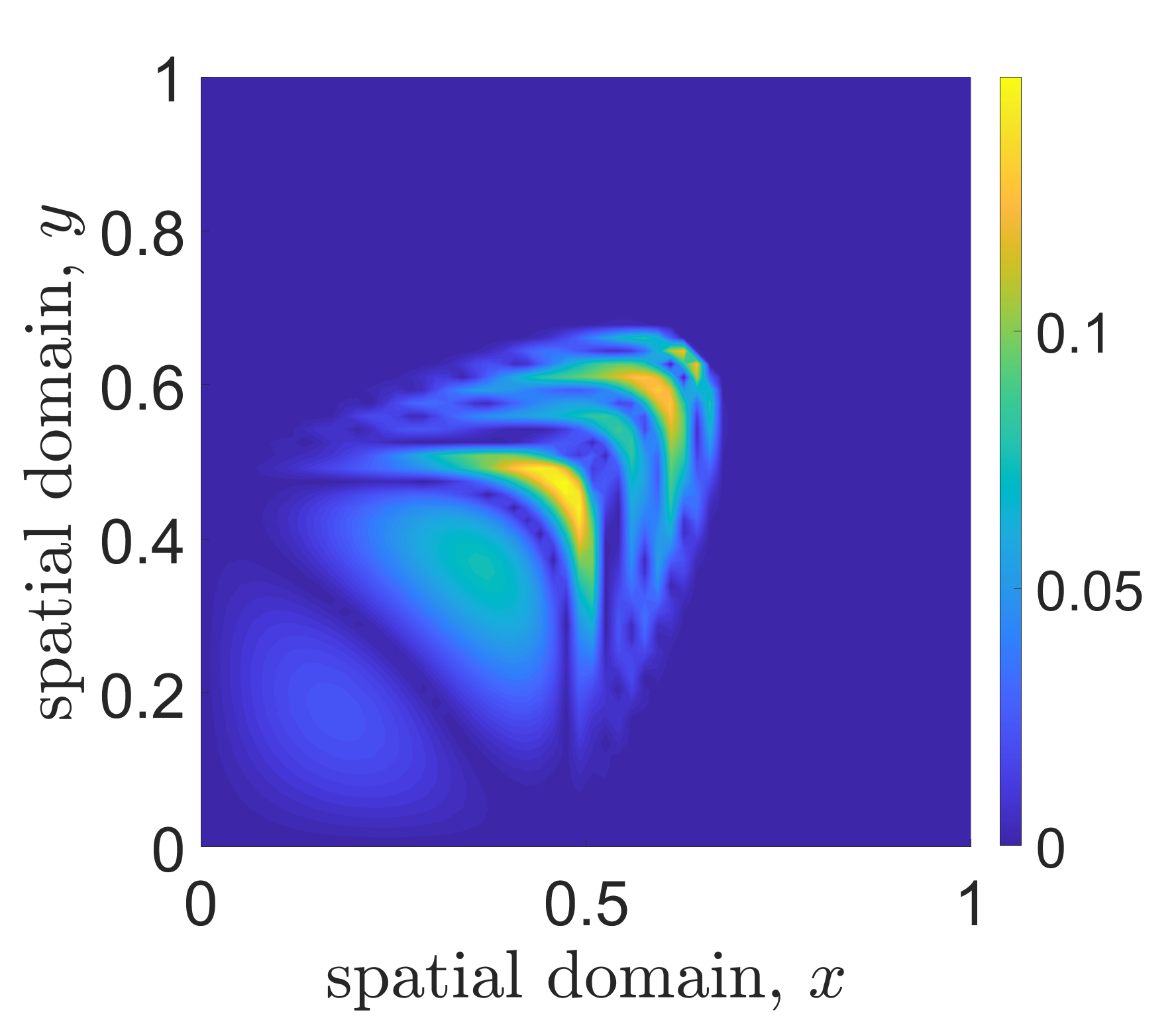}
}~~~~~~~ \hspace{-22pt}
\subfigure[BB-NN]{
  \includegraphics[width=0.25\textwidth]{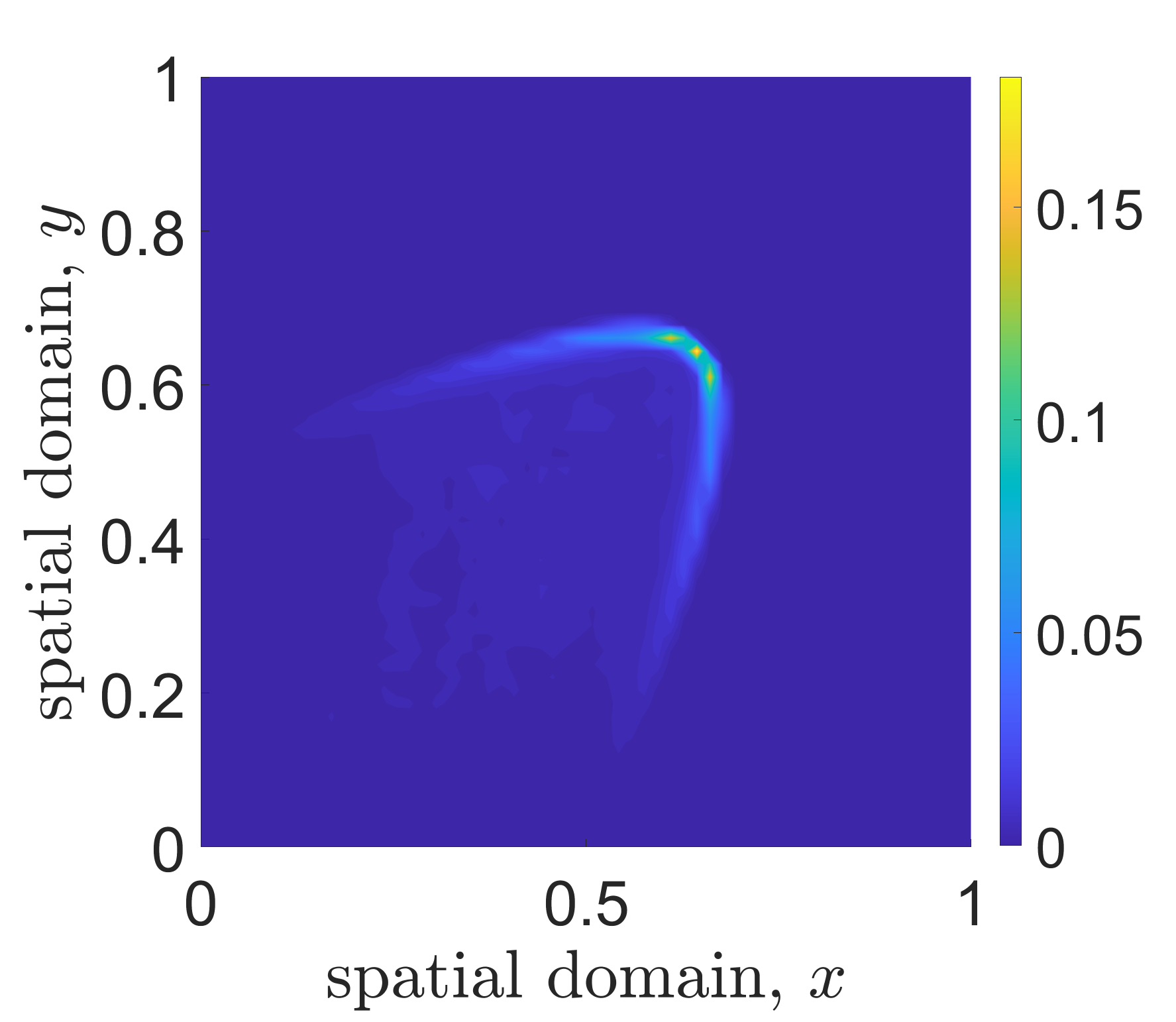}
}~~~~~~~
  \caption{(a) Solution snapshots, $u$, of FOM and absolute differences of
  (b) NM-LSPG-HR,
  (c) LS-LSPG-HR, and (d) BB-NN with respect to FOM solution at $t=2$.}
  \label{fg:2dburgersFOMROMS}
\end{figure}

We demonstrate the performance of our NM-ROMs (i.e., NM-LSPG and NM-LSPG-HR) by
comparing it with LS-ROMs (i.e., LS-LSPG and LS-LSPG-HR) that was first
introduced in \cite{carlberg2011efficient}. We solve the following parameterized
2D viscous Burgers' equation:
\begin{equation}\label{eq:2dburgers_eq} 
\begin{split}
\frac{\partial u}{\partial t} + u\frac{\partial u}{\partial x} + v\frac{\partial u}{\partial y} &= \frac{1}{Re}\left(\frac{\partial^2 u}{\partial x^2}+\frac{\partial^2 u}{\partial y^2}\right) \\ 
\frac{\partial v}{\partial t} + u\frac{\partial v}{\partial x} + v\frac{\partial v}{\partial y} &= \frac{1}{Re}\left(\frac{\partial^2 v}{\partial x^2}+\frac{\partial^2 v}{\partial y^2}\right) \\ 
(x,y) &\in \spaceDomain = [0,1] \times [0,1] \\
t &\in[0,2],
\end{split}
\end{equation}
with the boundary condition
\begin{equation}
    u(x,y,t;\mu)=v(x,y,t;\mu)=0 \quad \text{on} \quad \Gamma = \left\{(x,y)|x\in\{0,1\},y\in\{0,1\}\right\}
\end{equation}
and the initial condition
\begin{align}
    u(x,y,0;\mu) = \left\{
                                  \begin{array}{ll}
                                    \mu \sin{(2\pi x)}\cdot \sin{(2\pi y)} \quad &\text{if } (x,y) \in [0,0.5]\times [0,0.5] \\
                                    0 \quad &\text{otherwise}
                                  \end{array}
                                  \right.  \\
    v(x,y,0;\mu) = \left\{
                                  \begin{array}{ll}
                                    \mu \sin{(2\pi x)}\cdot \sin{(2\pi y)} \quad &\text{if } (x,y) \in [0,0.5]\times [0,0.5] \\
                                    0 \quad &\text{otherwise}
                                  \end{array}
                                  \right.
\end{align}
where $\mu\in \paramDomain=[0.9,1.1]$ is a parameter and $u(x,y,t;\mu)$ and
$v(x,y,t;\mu)$ denote the $x$ and $y$ directional velocities, respectively, with
$u:\spaceDomain \times [0,2] \times \paramDomain \rightarrow \RR{}$ and
$v:\spaceDomain \times [0,2] \times \paramDomain \rightarrow \RR{}$ defined as
the solutions to Eq.~\eqref{eq:2dburgers_eq}, and $Re$ is a Reynolds number
which is set $Re=10,000$. For this case, the FOM solution snapshot shows slowly
decaying singular values.  We observe that a sharp gradient, i.e., a shock,
appears in the FOM solution (e.g., see Fig.~\ref{fg:2dburgersFOMROMS}(a)).
We use $60\times 60$ uniform mesh with the backward difference scheme for the
first spatial derivative terms and the central difference scheme for the second spatial derivative terms. Then, we use the backward Euler scheme with time
step size $\Delta t = \frac{2}{\nt}$, where $\nt = 1,500$ is the number of time
steps.

For the training process, we collect solution snapshots associated with the
parameter $\mu \in \paramDomain_{train}=\{0.9,0.95,1.05,1.1\}$,
such that $\ntrain=4$, at which the FOM is solved. Then, the number of train
data points is $\ntrain\cdot(\nt+1)=6,004$ and $10\%$ of the train data are used
for validation purposes. We employ the Adam optimizer \cite{kingma2014adam} with
the SGD and the initial learning rate of $0.001$, which decreases by a factor of
$10$ when a training loss stagnates for $10$ successive epochs.  We set
the encoder and decoder hidden layer sizes to $6,728$ and $33,730$, respectively
and vary the dimension of the latent space from $5$ to $20$.  The weights and
bias of the autoencoder are initialized via Kaiming initialization
\cite{he2015delving}. The batch size is $240$ and the maximum number of
epochs is $10,000$. The training process is stopped if the loss on the
validation dataset stagnates for $200$ epochs. 

\begin{figure}[!htbp]                                                           
  \centering                                                                    
  \subfigure[Relative errors vs reduced
  dimensions]{                                              
  \includegraphics[width=0.48\textwidth]{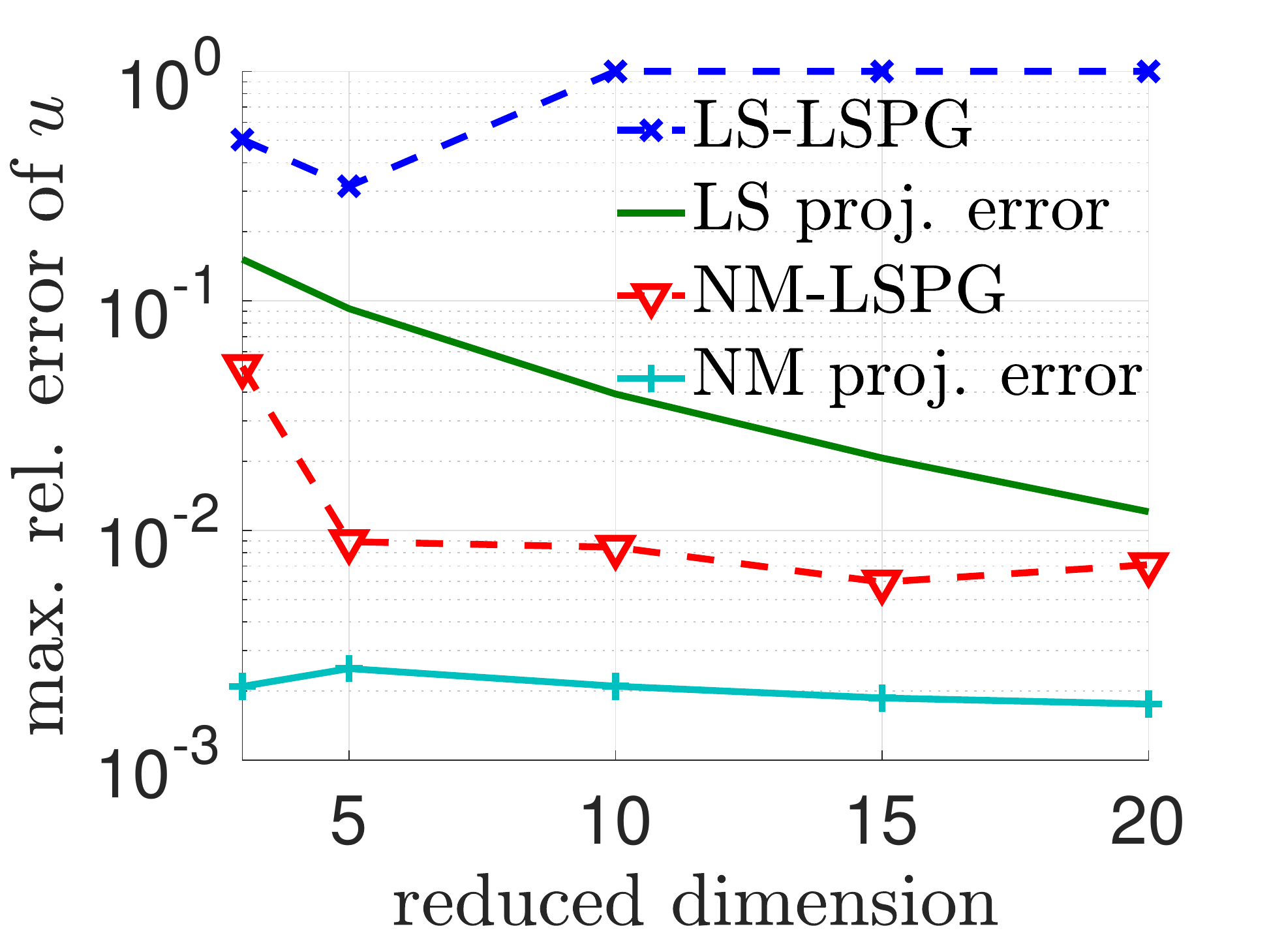}               
  }                                                                             
  \subfigure[Relative errors vs $\mu$]{                                              
  \includegraphics[width=0.48\textwidth]{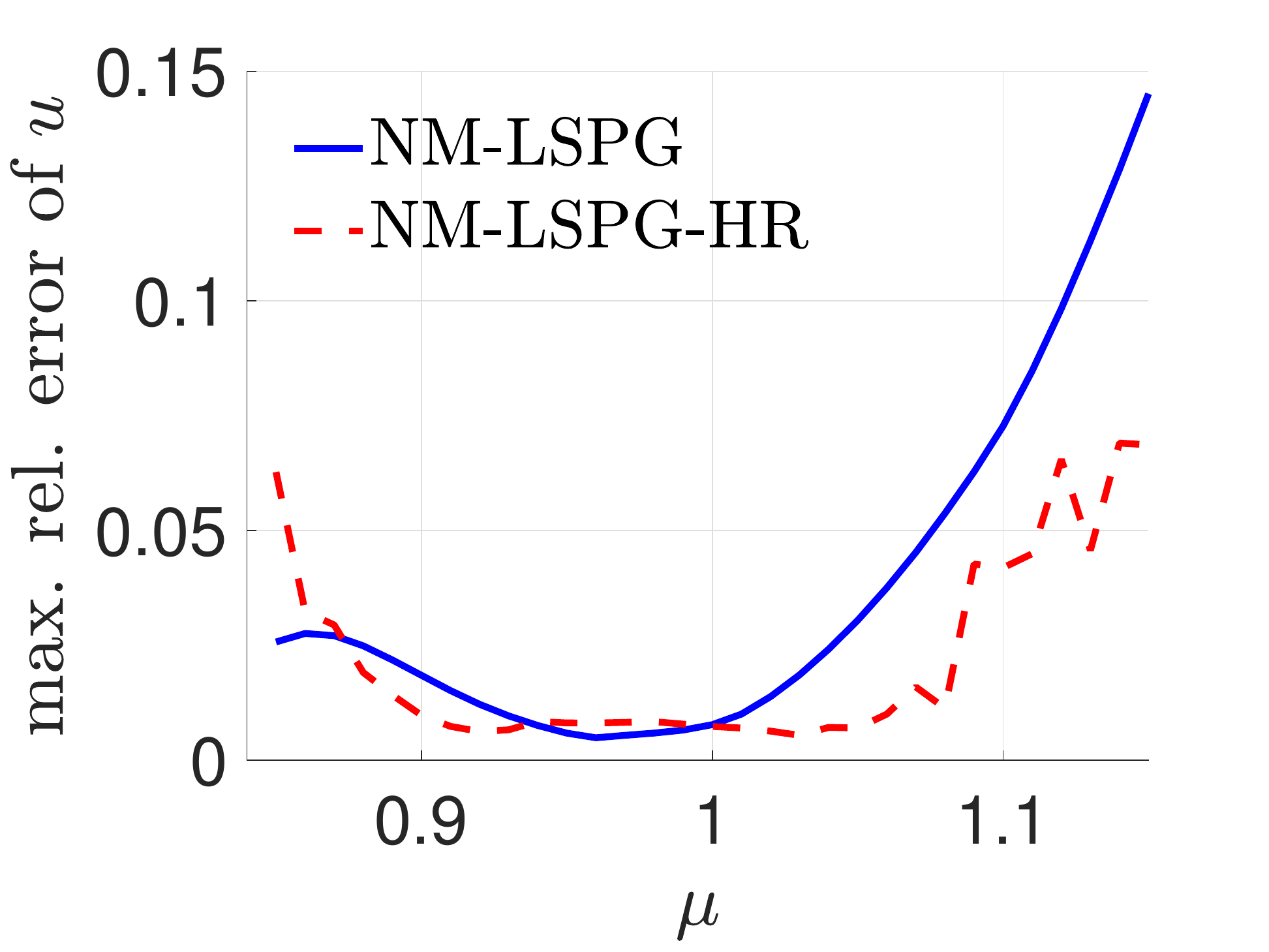}               
  }                                                                             
  \caption{The comparison of the NM-LSPG-HR and NM-LSPG on the maximum relative 
  errors. A maximum relative error that is $1$ means the ROM failed to solve the problem.}                                                               
  \label{fg:2dErrVSredDim}                                                   
\end{figure}

After the training is done, the NM-LSPG and LS-LSPG solve the
Eq.~\eqref{eq:2dburgers_eq} with the target parameter $\mu=1$, which is not
included in the train dataset. Fig.~\ref{fg:2dErrVSredDim} shows the relative
error versus the reduced dimension $\nbasisspace$ for both LS-LSPG and NM-LSPG.
It also shows the projection errors for LS-ROMs and NM-ROMs, which are the lower
bounds that any LS-ROMs and NM-ROMs can reach, respectively. As
expected the relative error for the NM-LSPG is lower than the one for the
LS-LSPG.  We even observe that the relative errors of NM-LSPG are even lower
than the LS projected error. 

We vary the number of residual basis and residual samples, with the fixed
number of training parameter instances $\ntrain=4$ and the        reduced
dimension $\nbasisspace=5$, and measure the wall-clock time. The results are
shown in Table~\ref{tb:2dDEIMtest}. Although the LS-LSPG-HR can achieve better
speed-up   than the NM-LSPG-HR, the relative error of the LS-LSPG-HR is too large
to be    reasonable, e.g., the relative errors of around $37 \%$. On the other
hand, the NM-LSPG-HR achieves much better accuracy, i.e., a relative error of
around $1   \%$, and a factor $11$ speed-up.

\begin{table}[!htbp]
\caption{The top 6 maximum relative errors and wall-clock times at different
  numbers of residual basis and samples which range from $40$ to
  $60$.}\label{tb:2dDEIMtest}
\centering
\resizebox{\textwidth}{!}{\begin{tabular}{|c|c|c|c|c|c|c|c|c|c|c|c|c|}
\hline
 & \multicolumn{6}{|c|}{NM-LSPG-HR} & \multicolumn{6}{|c|}{LS-LSPG-HR}\\
\hline
  Residual basis, $\nbasisres$ & 55 & 56 & 51 & 53 & 54 & 44 & 59 & 53 & 53 & 53 & 53 & 53\\
\hline
  Residual samples, $\nressample$ & 58 & 59 & 54 & 56 & 57 & 47 & 59 & 58 & 59 & 56 & 55 & 53\\
\hline
Max. rel. error (\%) & 0.93 & 0.94 & 0.95 & 0.97 & 0.97 & 0.98 & 34.38 & 37.73 & 37.84 & 37.95 & 37.96 & 37.97\\
\hline
Wall-clock time (sec) & 12.15 & 12.35 & 12.09 & 12.14 & 12.29 & 12.01 & 5.26 & 5.02 & 4.86 & 5.05 & 4.75 & 7.18\\
\hline
Speed-up & 11.58 & 11.39 & 11.63 & 11.58 & 11.44 & 11.71 & 26.76 & 28.02 & 28.95 & 27.83 & 29.61 & 19.58 \\
\hline
\end{tabular}}
\end{table}

Fig.~\ref{fg:2dburgersFOMROMS} shows FOM solutions at the last time and
absolute differences between FOM and other approaches, i.e., NM-LSPG-HR and LS-LSPG-HR with the reduced dimension being $\nbasisspace=5$ and a black-box NN
approach (BB-NN). The BB-NN approach is similar to the one described in \cite{DBLP:journals/corr/abs-1806-02071}. The main difference is that $L1$-norms and physics constraints were not used in our loss function. This approach gave a maximum relative error of $38.6\%$ and has a speed-up of $119$. While this
approach is appealing in that it does not require access to the PDE solver, the
errors are too large for our application. For NM-LSPG-HR, $55$
residual basis dimension and $58$ residual samples are used and for LS-LSPG-HR,
$59$ residual basis dimension and $59$ residual samples are used. Both FOM and
NM-LSPG-HR show good agreement in their solutions, while the LS-LSPG-HR is not
able to achieve a good accuracy.  In fact, the NM-LSPG-HR is able to achieve an
accuracy as good as the NM-LSPG for some combinations of the small number of
residual basis and residual samples.    

Finally, Fig.~\ref{fg:2dErrVSredDim}(b) shows the maximum relative error over
the test range of the parameter points.  Note that the NM-LSPG and NM-LSPG-HR
are the most accurate within the range of the training points, i.e., $[0.9,
1.1]$. As the parameter points go beyond the training parameter domain, the
accuracy of the NM-LSPG and NM-LSPG-HR start to deteriorate gradually. This
implies that  the NM-LSPG and NM-LSPG-HR have a trust region. Its trust region
should be  determined by each application. 

\section{Discussion \& conclusion}\label{sec:discussion-conclusion}
In this work, we have successfully developed an accurate and efficient NM-ROM.
We demonstrated that both the LS-ROM and BB-NN are not able to represent
advection-dominated or sharp gradient solutions of 2D viscous Burgers' equation
with a high Reynolds number. However, our new approach, NM-LSPG-HR, solves such
problem accurately and efficiently.  The speed-up of the NM-LSPG-HR is achieved
by choosing the shallow masked decoder as the nonlinear manifold and applying
the efficient hyper-reduction computation. Because the difference in the
computational cost of the FOM and NM-LSPG-HR increases as a function of the
number of mesh points, we expect more speed-up as the number of mesh points
becomes larger. 

Compared with the deep neural networks for computer vision and natural language
processing applications, our neural networks are shallow with a small number of
parameters. However, these networks were able to capture the variation in our 2D
Burgers' simulations. A main future work for transferring this work to more
complex simulations, will be to find the right balance between a shallow network
that is large enough to capture the data variance and yet small enough to run
faster than the FOM. Another future work will be to find an efficient way of
determining the proper size of the residual basis and the number of sample
points \textit{a priori}. To find the optimal size of residual basis and the
number of sample points for hyper-reduced ROMs, we relied on test results. This
issue is not just for NM-LSPG-HR, but also for LS-LSPG-HR.

\section*{Broader Impact}

The broader impact of this work will be to accelerate physics simulations to
improve design optimization and control problems, which require
thousands of simulation runs to learn an optimal design or viable control
strategy.  While this is not computationally feasible with high-fidelity FOMs,
the development of the NM-ROMs is an important step in this direction.

\begin{ack}
This work was performed at Lawrence Livermore National Laboratory and was
  supported by the LDRD program (project 20-FS-007).  Youngkyu was also
  supported for this work through generous funding from DTRA. Lawrence Livermore
  National   Laboratory is operated by Lawrence Livermore National Security,
  LLC, for the    U.S. Department of Energy, National Nuclear Security
  Administration under Contract DE-AC52-07NA27344 and LLNL-CONF-815209.  We
  declare that there were no conflicting interests of any type during the
  production of this research.
\end{ack}

\bibliographystyle{unsrt}
\bibliography{references}

\end{document}

%% file: commands.tex
\newcommand{\fontDiscrete}{\mathcal}

\newcommand{\bds}[1]{{\boldsymbol{#1}}}

\newcommand{\defeq}{:=}

\newcommand{\stateContinuousDummyEntryNo}{\boldsymbol{u}}

\newcommand{\stateContinuousDummy}[1]{\bds{\stateContinuousDummyEntryNo}} 
\newcommand{\stateTwoContinuousDummyEntryNo}{\boldsymbol{v}}

\newcommand{\stateTwoContinuousDummy}[1]{\bds{\stateTwoContinuousDummyEntryNo}}

\newcommand{\ndof}{\nspacedof}

\newcommand{\unitvec}{\bds{e}}
\newcommand{\unitvecArg}[1]{\unitvec_{#1}}

\newcommand{\argmin}[1]{\underset{#1}{\text{argmin}}}

\newcommand{\ndofPorts}[1]{n^p}

\newcounter{remctr}
\newcounter{propctr}
\newcounter{proposctr}
\newcommand{\RR}[1]{\ensuremath{\mathbb{R}^{ #1 }}}
\newcommand{\NN}{\mathbb{N}}
\newcommand{\RRplus}[1]{\ensuremath{\mathbb{R}_+^{ #1 }}}
\newcommand{\natNo}{\NN}
\newcommand{\nat}[1]{\natNo(#1)}

\newcommand{\BE}{\mathrm{BE}}

\newcommand{\identity}[1]{\boldsymbol I_{#1}}

\newcommand{\resSymb}{r}
\newcommand{\resRedSymb}{{\mathsf r}}
\newcommand{\fluxSymb}{f}
\newcommand{\paramSymb}{\mu}
\newcommand{\solSymb}{x}
\newcommand{\obliqueprojector}{\mathcal P}

\newcommand{\spaceSymb}{s}

\newcommand{\timeSymb}{t}

\newcommand{\weightmatSymb}{A}

\newcommand{\basismatspaceSymb}{\Phi}

\newcommand{\samplematSymb}{Z}
\newcommand{\samplevecSymb}{z}
\newcommand{\dummySymb}{v}

\newcommand{\nsmall}{n}

\newcommand{\nparam}{\nsmall_\mu}
\newcommand{\ntrain}{\nsmall_\text{train}}

\newcommand{\nressample}{\nsmall_{\samplevecSymb}}

\newcommand{\nbasisflux}{\nsmall_\fluxSymb}
\newcommand{\nbasisres}{\nsmall_\resSymb}

\newcommand{\nbasisspace}{{\nsmall_\spaceSymb}}

\newcommand{\nbig}{N}
\newcommand{\nspacedof}{\nbig_\spaceSymb}
\newcommand{\ntimedof}{{\nbig_\timeSymb}}

\newcommand{\paramDomain}{\mathcal D}

\newcommand{\totaltime}{T}

\newcommand{\spatialSubspace}{\fontDiscrete S}

\newcommand{\res}{\boldsymbol \resSymb}

\newcommand{\resArg}[1]{\boldsymbol \resSymb^{#1}}

\newcommand{\resRedApprox}{\tilde{\boldsymbol \resRedSymb}}
\newcommand{\resRedApproxArg}[1]{\resRedApprox^{#1}}

\newcommand{\resn}{\resArg{n}}

\newcommand{\param}{\boldsymbol \paramSymb}

\newcommand{\paramDummy}{\boldsymbol \nu}
\newcommand{\dt}{\Delta \timeSymb}

\newcommand{\sol}{\boldsymbol \solSymb}

\newcommand{\flux}{\boldsymbol \fluxSymb}
\newcommand{\fluxArg}[1]{\flux_{#1}}

\newcommand{\solDummy}{\boldsymbol w}

\newcommand{\timeArg}[1]{t^{#1}}
\newcommand{\timeDummy}{\tau}
\newcommand{\timeDomain}{[0,\totaltime]}

\newcommand{\solArg}[1]{\sol_{#1}}
\newcommand{\solFunc}{\sol}

\newcommand{\solFuncArg}[1]{\solFunc(t^{#1};\param)}

\newcommand{\solapproxFunc}{\solapprox}
\newcommand{\solapproxFuncArg}[1]{\solapproxFunc(t^{#1};\param)}

\newcommand{\redsolapproxFunc}{\redsolapprox}
\newcommand{\redsolapproxFuncArg}[1]{\redsolapproxFunc(t^{#1};\param)}

\newcommand{\testbasisArg}[2]{\boldsymbol \xi_{ij}}

\newcommand{\solapproxArg}[1]{\solapprox_{#1}}

\newcommand{\solapprox}{\tilde\sol}

\newcommand{\redsolapprox}{\hat\sol}

\newcommand{\redsolapproxArg}[1]{\redsolapprox_{#1}}

\newcommand{\dummy}{\boldsymbol {\dummySymb}}
\newcommand{\reddummy}{\hat{\dummy}}

\newcommand{\redres}{\hat{\res}}

\newcommand{\basismatspace}{\boldsymbol{\basismatspaceSymb}}

\newcommand{\basismatres}{\basismatspace_\resSymb}

\newcommand{\basisresvecArg}[1]{\boldsymbol \phi_{r,#1}}

\newcommand{\weightmat}{\boldsymbol{\weightmatSymb}}

\newcommand{\samplematNT}{\boldsymbol \samplematSymb}
\newcommand{\samplemat}{\samplematNT^T}

\newcommand{\unscaledEncoder}{\boldsymbol{E}}
\newcommand{\unscaledDecoder}{\boldsymbol{D}}

\newcommand{\scaledDecoder}{\boldsymbol{g}}

\newcommand{\spaceDomain}{\Omega}

\newcommand{\nt}{n_t}

\newcommand{\RN}[1]{%
  \textup{\uppercase\expandafter{\romannumeral#1}}%
}
